\documentclass[11pt]{article}
\usepackage{amssymb}
\usepackage{}
\usepackage{pifont}
\usepackage{amscd}
\usepackage{amsfonts}
\usepackage{amsmath, amsthm}
\usepackage{amssymb, amsxtra}

\usepackage{graphicx}

\usepackage{multicol}
\usepackage{epstopdf}
\usepackage{epsf}
\usepackage{float}
\usepackage{extarrows}
\usepackage[a4paper, text={.8\paperwidth,.83\paperheight}, ratio=1:1]{geometry}
\usepackage[format=hang,font=small,textfont=it]{caption}
\usepackage{authblk}

\usepackage[colorlinks=true, linkcolor=blue, citecolor=blue]{hyperref}
\usepackage[english]{babel}
\usepackage{latexsym}

\usepackage{tikz}
\usetikzlibrary{calc}
\usepackage[tight]{subfigure}
\usepackage{mathrsfs}
\usepackage{color}
\definecolor{greenbean}{RGB}{199,237,204}
\newtheorem{thm}{Theorem}[]

\newtheorem{lemma}[thm]{Lemma}

\newtheorem{eg}{Example}[]

\newtheorem{Def}[thm]{Definition}

\def \ta{\tau}

\def \ta1{\tau_1}

\begin{document}
\title{A note on a question of Shioda about integral sections
\footnotetext{\hspace{-1.8em} Email address: Jia-Li Mo:  mojiali0722@126.com;\\
2020 Mathematics Subject Classification.  11G05, 14H10, 14J27.}}

\author[1]{Jia-Li Mo}
\affil[1]{\small{Department of Mathematics, Soochow University, Shizi RD 1, Suzhou 215006, Jiangsu, China}}

\date{}
\maketitle

\abstract{We consider a rational elliptic surface with a relatively minimal fibration. We compute the number of integral sections in the above rational elliptic surface. As an application, we obtain an estimate of polynomial solutions of some equations.}


\section{Introduction}\label{outline}

Let $C$ be a smooth curve over a filed $k$ (an algebraically closed field of characteristic zero), and $K$ be its function field.
Now, we give an elliptic curve over $K$,
a $K$-rational point $P$  on the original elliptic curve over $K$ is called a section.

In this paper,  we  use the language of algebraic surface.  Let $S$ be a smooth projective surface having a relatively minimal elliptic fibration $f : S  \rightarrow C$ with the zero-section $O$ over a curve $C$. Assume that $S$ has at least one singular fibre, then the group $E(K)$ of $K$-rational points is finitely generated (Mordell-Weil Theorem). It can be identified with the group of sections of $f$. For each $P$ in $E(K)$, we denote by $(P)$ the image curve of the corresponding section $C\rightarrow S$, the curve $(P)$ may be also called a {\bf section} without confusion.
An element $P$ of $E(K)$ is called  integral if $(P)$ is disjoint from the zero-section $(O)$. Let $\mathcal{P}$ be the set of all integral sections:
$$ \mathcal{P}=\{P\in E(K)| (P)\cap (O)=\varnothing \}$$

The easy analogue of Siegel's theorem (\cite{Sie29}) about the finiteness of integral sections of elliptic curves can refer to \cite{S2}: {\bf Any elliptic surface admits only finitely many integral sections}.
%For a rational elliptic surface, there is   much more precise result about the number of integral sections (\cite{TGM1}): {\bf For any rational elliptic surface, there are at most 240 integral sections.}

Moreover, Shioda gave an interesting open question in \cite{TGM1}: {\bf Given an elliptic surface $S$
over $\mathbb{P}^1$ of arithmetic genus $\chi(S)$, with the generic fibre
$F$, what is the number of integral sections? }

Shioda in  \cite{TGM1} proved that the number is at most 240  for rational surfaces, and he computed the exact number for $r=8,7,6$. It is well known that, this open problem involves not only combinations but also the arithmetic structures of elliptic surface if $r\leq 5$.
In our paper,  we totally solve this open question for rational surfaces. We give the exact number for integral sections of every case. In fact, we consider the following sets $\mathcal{P}_i=\{P_i\in E(K)| (P_i)\cap (O)=\varnothing~{\rm and}~<P_i, P_i>=k_i \}$ in our paper. Now, we state the main theorem:

\begin{thm}\label{thm1}
The structures of integral sections of every rational elliptic surface  are given in table \ref{tab1:my_label}--\ref{tab10:my_label}.
In these tables, the rows are arranged in the following order and  we also give other indices, they descend with
respect to:
\end{thm}

\begin{itemize}
  \item (1)~$r$ ----- the rank~ of~~ $E(K)$; (2)~ $T$ ----- the trivial lattice; (3)~$v(T)$ ----- the number of roots in the trivial lattice $T$;
   (4)~$<P_i,P_i>$ ----- the height pairing; (5)~$m(P_i)$ ----- the number of  combinatorial multiplicity of $P$; (6)~${\mathrm \#} {\mathcal{P}_i}$ ----- the number of integral sections of different  height pairing.
\end{itemize}

%As we know, the rank problem for arbitrary elliptic curve over $\mathbb{C}(t)$ is still open. It is not easy to construct a elliptic curve with large rank (see \cite{Ul}).
%In general, an elliptic surface has more sections, then it will have more integral  sections.
%By the deformation theory of surfaces, Shioda proved that: let $f_i: S_i\rightarrow \mathbb{P}^1$ be two  elliptic fibrations. If $MW(f_1)>MW(f_2)\geq6$, then  $MW(f_1)>MW(f_2)$.
%From Theorem \cite{thm1}, we can shaper our result to 4.  And 4 is the sharper.
Estimate the number of integer solutions of an  equation is an interesting question.
Similarly, we consider  the number of polynomial solutions of an equation.  As an application of Theorem \ref{thm1}, we try to use discriminant to estimate the number of  polynomial solutions  for the following Weierstrass equations:

\begin{thm}\label{thm2}
For $$ y^2= x^3+a(t)x^2+b(t)x+c(t),$$
where $a(x), b(x), c(x)$ are polynomials. And deg~$(a(t))\leq2$, deg~$(b(t))\leq4$, deg~$(c(t))\leq6$.
Let the number of  polynomial solutions  be $n$, and the discriminant be $\Delta(t)$. Then
\begin{itemize}
  \item $0\leq n \leq 3$, if $\Delta(t)$ has only one solution, which is not 12-fold root.
  \item $0\leq n\leq N_0$ (where $N_0\geq 26$), if $\Delta(t)$  satisfies one of the following conditions:
  \begin{itemize}
    \item (i)~it has only two different solutions.
    \item (ii)~it has only three different solutions and deg~$(\Delta(t))=12$.
  \end{itemize}
  \item $n \geq 1$, if $\Delta(t)$  satisfies one of the following conditions:
  \begin{itemize}
    \item (i)~it has  more than three different solutions.
    \item (ii)~it has only three different solutions and deg~$(\Delta(t))<12$.
    \item (iii)~it has $3\leq \mathrm{deg}~\Delta(t) \leq 9$.
  \end{itemize}
\end{itemize}
\end{thm}

\paragraph{Acknowledgements:}
I thank Dr. Cheng  Gong for useful discussions about the elliptic fibrations. I also thank the reviewers for their good advise. This paper was  partly supported by the Natural Science Foundation of Jiangsu Province~(BK 20211305), and Postgraduate Research \& Practice Innovation Program of Jiangsu Province~(KYCX22\_3180).

\section{Definitions and terminology}\label{method}

By a rational elliptic surface, we mean a (smooth projective) rational surface over
$k$, say $S$, which is given with a relatively minimal elliptic fibration $f : S \rightarrow C$.
Now, we study the number of integral sections  of rational elliptic surfaces with sections.

First, we recall and fix some standard notation in dealing with Mordell-Weil groups (see \cite{Shio91}). The reader can refer to \cite{OS, S2} for more details.
\begin{itemize}
  % \item  $(P)$ : the curve on $S$ determined by a section $P$, esp. $(O)$  is the zero-section viewed as a curve on $S$.
%    \item  $(P,O)$ : the intersection number of $(P)$ and $(O)$.
% \item $E(K)$ : the Mordell-Weil group, i.e. the group of $K-$ rational points of a elliptic curve $E$ over $K$ with the zero-section $O$, which is finitely generated.
\item  $<P, Q>$ : the height pairing $(P, Q \in E(K))$, as defined by Shioda in \cite{S2}.
 This is a symmetric bilinear pairing which is positive-definite modulo torsion.

  \item  $R : = \{v \in C| f^{-1}(v)~ is~ reducible \}$.
  \item  $m_v$ : the number of irreducible components $\Theta_{v,i}$ of fiber $F_v$ of an elliptic surface.
  \item  $\Theta_{v,i}$~$(0\leq i\leq m_{v-1})$ : irreducible components of $f^{-1}(v)(v\in R)$, with $i=0$ corresponding to the identity component, which intersects with zero section $(O)$.
  \item  $T_v$ : the lattice generated by $\Theta_{v,i}~(i>0)$ with the sign changed, this is  a root lattice of type $A,D,E$ determined by the type of reducible fiber $f^{-1}(v)$.
 \item  $T=\bigoplus_{v\in R}T_v$ : the trivial lattice.
  \item  $E(K)^{0}$ : the narrow Mordell-Weil lattice of $E/K$, and $ E(K)^{0}=\{ P\in E(K)|(P)~{\rm meet}~\Theta_{v,0},\\~{\rm for~all}~v \in R \}$.  This is a certain subgroup of finite index in $E(K)$.
\item  ${E(K)^{0}}^{\vee}$ : the dual lattice of $E(K)^{0}$(= free part in $E(K)$).

    \item $v(T)$ : the number of roots in the trivial lattice $T$, especially, if we fix $T$, then value of $v(T)$ can be computed by {\cite[Table ~2.2]{SS}}.

% \item $contr_{v}(P,Q)$: Let $P,Q\in E(K)$. Then we define the local contribution from the singular fiber at $v \in R$. and $contr_{v}(P)=contr_{v}(P,P)$.

% \item $m(P)$: the combinatorial multiplicity $m(P)$ which is defined as the number of the distinguished roots in the root graph associated with $P$
% {\cite[Theorem 3.4]{TGM1}}
\end{itemize}
% $dim_kR/I$ : the linear dimension $dim_kR/I$ is equal to the number of integral sections on $S$ counted with multiplicities. Assume that $S$ is a rational elliptic surface, The linear dimension $dim_k R/I$ is equal to
%      $240-v(T)$, where $v(T)$ is the number of roots in the trivial lattice $T$.

We know that the Mordell-Weil group $E(K)$ of an elliptic surface
$S$ is canonically isomorphic to the quotient group $NS(S)/ Triv(S)$, where $Triv(S)=T\bigoplus<O,F>$. It has a height pairing  on $E(K)$.

For any $P, Q \in E(K)$, let
$<P, Q> = (\varphi(P). \varphi(Q))$.
Then it defines a $\mathbb{Q}$-valued symmetric bilinear pairing on $E(K)$ which induces the
structure of a positive-definite lattice on $E(K)/E(K)_{tors}$. The lattice $(E(K)/E(K)_{tors}, <\cdot,\cdot>)$ will be called the
Mordell-Weil lattice of the elliptic curve $E/K$, or of the elliptic surface $f : S \rightarrow C$
over $k$.% We will often abbreviate the Mordell¨CWeil lattice as $MWL(S)$
%or even $MWL$ if no confusion is likely.

The following theorem gives us  the explicit formula (see \cite{S2}) for the height pairing.
\begin{thm} {\bf (Shioda,1990)}\label{thm3}
For any $P, Q \in E(K)$, we
have
$$<P, Q> = \chi(S) + (P.O) + (Q.O)-(P.Q)-\sum_{v\in {R}}{contr_{v}(P, Q)}, $$
in particular,
$$<P, P>= 2\chi(S) + 2(P.O)-\sum_{v\in {R}}{contr_{v}(P)}. $$

\end{thm}

Now, suppose that $(P)$ intersects $\Theta_{v,i}$ and $(Q)$ intersects $\Theta_{v,j}$ and assume that
$i \geq 1, j \geq 1$. Then the contribution terms $contr_{v}(P, Q)$ and $contr_{v}(P)$ are given in  {\cite[Table (8.16)]{S2}}.

Since we also give another index  in Theorem \ref{thm1}, we're going to briefly introduce $m(P)$.
It is the combinatorial multiplicity which is defined by Shioda in \cite{TGM1}.

\begin{Def}
For any $P \in E(K)$, let $R(P)$ denote the subset of $v \in R$ such that $(P)$ intersects
some non-identity component $\Theta_{v,i}(i\neq0)$ of $f^{-1}(v)$.
The root graph associated with $P$, denoted by $\Delta(P)$,
is the connected graph with the vertices
$D(P),\Theta_{v,i}(v\in R(P),i\neq0)$,  where $D(P):=(P)-(O)-F$
and two vertices $\alpha,\beta$ are connected by an edge iff
the intersection number $\alpha\cdot\beta=1$. By a distinguished
root of $\Delta(P)$, we mean a linear combination of the
vertices of the form:
$D_p=D(P)+\sum_{v,i}n_{v,i}\Theta_{v,i} (n_{v,i}\in \mathbb{Z},\geq 0)$
satisfying ${D_p}^2=-2$.
Further we denote by $m(P)$
the number of distinguished roots in the root graph
$\Delta(P)$, and call it the combinatorial multiplicity of $P$.
\end{Def}

Moreover, Shioda gave the following theorem for an rational elliptic surfaces over $\mathbb{P}^1$ in \cite{TGM1}:
\begin{thm}\label{theorem5}
$$\sum_{P\in {\mathcal{P}}}m(P)
=240-v(T)$$

\end{thm}

% the third row gives
%the value of ${contr_{v}(P, Q)}$  for the case $i = j$ as well as ${contr_{v}(P)}$ , and the fourth
%row concerns the case $i < j$ (interchange the order of $P, Q$ if necessary).
%(About contr can refer book table6.1)

%\item  $f : S \rightarrow C=\mathbb{P}^1$: the associated elliptic surface (the Kodaira-Neron model) of $E/K$ with at least one singular fiber. A $K-$ rational point $P\in E(K)$ is identified with a section of $f$.

%One of the obvious advantages of this result is that when computing the number of integral section on the Mordell-Weil group it suffices to use the following proposition.

\section{Proof of the Theorem 1}\label{proof}

In this section, we will give the full proof of  Theorem \ref{thm1}. First, for our purposes, we take into account a lemma as following:

\begin{lemma}\label{lemma1}
On a rational elliptic surface $S$, $P$ is a section and its height pairing $<P, P>=2$. Then $P$ is an integral section if and
only if $P\in E(K)^0$.
\end{lemma}

{\bf Proof:} If $P$ is an integral section,  $(O.P)=0$ by definition. then we get $\sum_{v\in {R}}{contr_{v}(P)}=0$ by the height pairing formula,
that is $P\in E(K)^0$.

If $P\in E(K)^0$, we have $\sum_{v\in {R}}{contr_{v}(P)}=0$, from height pairing $2=<P, P>= 2 + 2(P.O)-0$, it follows that $(O.P)=0$, then $P$ is an integral section. $\Box$

Here, we give another lemma  which  is used later:
\begin{lemma}\label{lemma2}
Let $P$ be a section of the Mordell-Weil group $E(K)$,  $P\in E(K)$ and $P \notin E(K)^0$ (we use $P\in E(K)-E(K)^0$ later). There exists $P'\in {E(K)^{0}}^{\vee}$ such that
$<P, P>=<P', P'>$. Moreover, ${\mathrm \#} {\mathcal{P}}\leq{\mathrm \#} P' \cdot \mid {E(K)_{tor}}\mid$. Where  ${\mathrm \#} {\mathcal{P}}$ means the number of integral sections with the same height pairing value $<P,P>$, and ${\mathrm \# P' }$ means the number of  sections in  ${E(K)^{0}}^{\vee}$ with the same height pairing value $<P,P>$.
\end{lemma}

{\bf Proof:}  $P\not\in {E(K)^{0}}^{\vee}$, there exists a torsion element $\alpha \in {E(K)_{tor}}$ of order $n$, such that $P-\alpha \in {E(K)^{0}}^{\vee}$. let $P'=P-\alpha$, then we have
$<P', P'>=<P-\alpha, P-\alpha>=<P,P>$. So it gives rise to ${\mathrm \#} P\leq{\mathrm \#} P' \cdot \mid {E(K)_{tor}}\mid$. $\Box$

%\begin{lemma}\label{lemma3}
%Under the condition that the Mordell-Weil group of rational elliptic surface is $E(K)\cong E(K)/E(K)_{tor}\bigoplus \mathbb{Z}_2$, a integral section $P\in E(K)$, Then ${\mathrm \#} P= {\mathrm \#} (P+\bar{1})$.
%\end{lemma}
%
%Proof: First, if a section $P\in E(K)/E(K)_{tor}$ pass through a fixed component of every reducible fibre $f^{-1}(v)$. We denote this intersection ways by $M$, it is easy to verify that $P+ \bar{1}$  must pass through a different component of every reducible fibre than $P$, Where $\bar{M}$ stands for intersection ways of $P'$ and every reducible fibre. In what follows, we construct two maps: $M\rightarrow \bar{M}$ and $\bar{M}\rightarrow M$, then $P \in M \Rightarrow P+ \bar{1}\in \bar{M}$, and $P \in \bar{M} \Rightarrow P + \bar{1}\in M$. Therefore, $M\cong \bar{M}$, it gives rise to  ${\mathrm \#} P= {\mathrm \#} (P+\bar{1})$.

Now, we start to give the full proof of Theorem 1---compute
the combinatorial multiplicity of integral sections and the number of integral sections of every rational elliptic surface.
We know that the Mordell-Weil groups $E(K)$ of rational surfaces have 74 types (see \cite{OS}). First, we define two situations as follows:

\begin{itemize}
  \item Situation ($S1$): There exists two kinds of $\sum_{v\in {R}}{contr_{v}(P)}$. The difference between the two numbers is $2$.
  \item Situation ($S2$): For sections $P_1, P_2$, they have the same height pairing value which corresponds to a different combinatorial multiplicity of $P_i$. In fact, the intersection ways of $P_i$ with the fibre are also different.
\end{itemize}

According to the above two situations, we shall divide the proof into four cases:

{\bf Case $\uppercase\expandafter{\romannumeral 1}$:} ~1--20, ~23, ~25--27, ~30--33, ~35--37, ~39, ~40, ~43--47,
~49, ~50, ~55, ~56, ~62--65, ~67--69 ~and~ 73 have no Situations ($S1$) and $(S2)$.

{\bf Case $\uppercase\expandafter{\romannumeral 2}$:}  ~24, ~41, ~42 ~and ~61  have only Situation $(S1)$.

{\bf Case $\uppercase\expandafter{\romannumeral 3}$:}  ~21, ~22, ~28, ~29, ~34, ~48, ~51, ~52, ~54, ~57, ~58, ~66, ~70--72 ~and ~74  have only Situation $(S2)$.

{\bf Case $\uppercase\expandafter{\romannumeral 4}$:}  ~38, ~53, ~59 ~and ~60 have both Situations $(S1)$ and $(S2)$.

{\bf The proof of  Case $\uppercase\expandafter{\romannumeral 1}$:}
In 74 types of $E(K)$, the most basic case for rational elliptic surfaces, the case $rank (T) \leq 2$ which
corresponds to the case of higher Mordell-Weil rank $r$ = 8, 7 or 6~(cases 1-4). They were verified in \cite{TGM2}.

\begin{figure}[ht]
\begin{center}
\scalebox{0.8}{\includegraphics{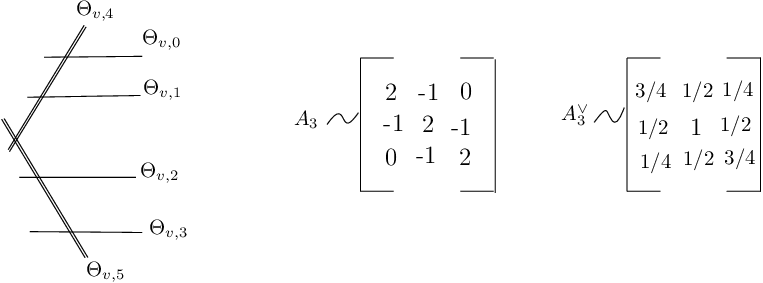}}
\end{center} \caption{}\label{No16}
\end{figure}

For the remaining cases, we need to choose another way to accomplish them and make use of Shioda's method on the Mordell-Weil groups.
With the given concepts, we are now ready to prove the case No. 16, ~$T=D_5$, where $v(T)=40$ and Mordell-Weil group $E(K) \cong A^{\vee}_3$.  It is a representation  for case $\uppercase\expandafter{\romannumeral 1}$. For a reducible fibre is ${I_1}^*$ (see Figure \ref{No16}). First, according to different the height pairing value of sections, we compute $m(P_i)$ for  integral section respectively.

By the height formula in Theorem \ref{thm3},
we have for any $P \in E(K), <P, P>= 2\chi(S) + 2(P.O)-\sum_{v\in R}{contr_{v}(P)}$, then we find the height pairing:
$$<P, P>= 2 + 2(P.O)-\left\{
\begin{array}{ll}
0, & \hbox{$P\in E(K)^0$, $P$ intersects $\Theta_{v,0}$;} \\
1, & \hbox{$P\in E(K)-E(K)^0$, $P$ intersects $\Theta_{v,1}$;} \\
5/4, & \hbox{$P\in E(K)-E(K)^0$, $P$ intersects $\Theta_{v,2}$;} \\
5/4, & \hbox{$P\in E(K)-E(K)^0$, $P$ intersects $\Theta_{v,3}$.}
\end{array}
\right.
$$
If we suppose that $ P_i \in \mathcal{P}_i ~(i=1,~2,~3, 4)$ are integral sections, since intersection ways of $P_3$ and $P_4$ with the fibres are symmetric, we can merge  $P_3$ and $P_4$ into one new set.  Then there are sections of height $<P_i, P_i>=2, ~1, ~3/4, ~(i=1,~2,~3)$. By inspection of the height formula, this requires that these sections are disjoint from the zero section and meet ${I_1}^*.$

First, by a distinguished root of the root graph associated with $P_i$, we mean a linear combination of the
vertices of the form: $D_p=D(P)+\sum_{v, i}{n_{v,i}\Theta_{v,i}}~(v\in R(P), ~i\neq 0)$ and $(n_{v,i}\in \mathbb{Z}, \geq 0)$ satisfying ${D_p}^2=-2$.
We have $((P)-(O)-(F)+\sum_{v, i}{n_{v,i}\Theta_{v,i}})^2=-2 ~~(i=1, \cdots, 5),$
~$F^2=0$  by \cite{TGM2}. For any $P \in E(K),$ we have $P^2=-\chi(S)=-1$.
Since it just has one singular fibre ${I_1}^*$, which is a reducible singular fibre. Then every component $\Theta_{v,i}$ of $F_v$
is a smooth rational curve which has self-intersection number $(\Theta_{v,i})^2 = -2,$ and  $(\Theta_{v,i}.\Theta_{v,j})=1 ~((i,j)=(1,4),~(2,5),~(3,5),~(4,5))$ which can be read off from the above Figure \ref{No16}.

 When $<P_1, P_1>=2$, we obtain
\begin{equation}\label{equation1}
\sum\limits_{i=1}^5 (n_{v,i})^2-(n_{v,1}n_{v,4})-(n_{v,2}n_{v,5})-(n_{v,3}n_{v,5})-(n_{v,4}n_{v,5})=0.
\end{equation}
It just has only one solution: $(0,0,0,0,0)$. Thus $m(P_1)=1$.

When $<P_2, P_2>=1,$ in the same way of above. Since  $(P_2 \Theta_{v,1})= 1$, we get equation
\begin{equation}\label{equation11}
\sum\limits_{i=1}^5 (n_{v,i})^2-(n_{v,1}n_{v,4})-(n_{v,2}n_{v,5})-(n_{v,3}n_{v,5})-(n_{v,4}n_{v,5})-n_{v,1}=0.
\end{equation}
This equation can be computed the number of solutions, which has 10 integer solutions:
\begin{itemize}
  \item  $(n_{v,1}, n_{v,2}, n_{v,3}, n_{v,4}, n_{v,5})$:
  \item (2,1,1,2,2), (1,1,1,2,2), (1,1,1,1,1), (1,1,1,1,2), (1,1,0,1,1)
  \item (1,0,1,1,1), (1,0,0,1,0), (1,0,0,1,1), (1,0,0,0,0), (0,0,0,0,0)
\end{itemize}
It means  $m(P_2)=10$.

When $<P_3, P_3>=3/4,$ because intersection number of $P_3$ and $\Theta_{v,2}$ or  $\Theta_{v,3}$ is 1, we have
\begin{equation}\label{equation111}
\sum\limits_{i=1}^5 (n_{v,i})^2-(n_{v,1}.n_{v,4})-(n_{v,2}.n_{v,5})-(n_{v,3}.n_{v,5})-(n_{v,4}.n_{v,5})-n_{v,3}=0.
\end{equation}
$m(P_3)=16$,  two intersection ways have same $m(P_3)$ and its integer solutions as follows:

\begin{itemize}
  \item (1,2,2,2,3), (1,1,2,2,2), (1,1,2,2,3), (1,1,2,1,2), (1,1,1,2,2), (1,1,1,1,1)
  \item (1,1,1,1,2), (1,0,1,1,1), (0,1,2,1,2), (0,1,1,1,1), (0,1,1,1,2), (0,1,1,0,1)
  \item (0,0,1,1,1), (0,0,1,0,0), (0,0,1,0,1), (0,0,0,0,0)
\end{itemize}

In the following, we compute the number of their integral sections.

When $<P_1, P_1>=2$, with the notation of Lemma \ref{lemma1}, it has $P_1\in E(K)^0 \cong A_3$, the Gram matrix of $A_3$ as above Figure 1. We know that $rank(E(K))= rank(E(K)^0)=3$, assume that $e_1, e_2, e_3$ are basis of $E(K)^0$ such that  $P_1=ae_1+be_2+ce_3 ~(a, b, c \in \mathbb{Z})$. According to $<P_1, P_1>=2$, it follows that $2a^2-2a b+2b^2-2b c+2c^2=2$, which has 12 integer solutions $(a,b,c)$ below by computing the equation. It means ${\mathrm \#} {\mathcal{P}_1}=12$, this is easily verified that these 12 integral solutions are exactly 12 integral sections.
\begin{itemize}
  \item (0,0,$\pm1$), (0,$\pm1$,0), ($\pm1$,0,0), (0,1,1), (0,-1,-1), (1,1,0), (-1,-1,0), (1,1,1), (-1,-1,-1)
\end{itemize}

When $<P_2, P_2>=1, P_2\in E(K)-E(K)^0$, we assume that $e'_1, e'_2, e'_3$ are  basis of $E(K)$ such that  $P_2=ae'_1+be'_2+ce'_3~ (a, b, c \in \mathbb{Z})$. According to $<P_2, P_2>=1$, it follows that $3/4a^2+ab+b^2+1/2ac+bc+3/4c^2=1$, which has 6 integer solutions. We need to check that whether there indeed are 6 integer solutions corresponding to 6 integral sections by an explicit height pairing. We know $\sum_{v\in {R}}{contr_{v}(P_2)}=2\chi(S) +2(O.P_2)-<P_2,P_2>$ and $max(\sum_{v\in {R}}{contr_{v}(P_2))=5/4},$ so $2+2(O.P_2)-1\leq 5/4$, moveover $(O.P_2)\leq 1/8$. We get $(O.P_2)=0$. In fact, we have 6 integral sections  as follows:

\begin{itemize}
  \item (0,$\pm1$,0), (1,0,-1), (-1,0,1), (-1,1,-1), (1,-1,1)
\end{itemize}

When $<P_3, P_3>=3/4$, $e'_1, ~e'_2, ~e'_3$ are basis of $E(K)$ such that  $P_2=ae'_1+be'_2+ce'_3 (a, b, c \in \mathbb{Z})$. According to $<P_3, P_3>=3/4$, it follows that $3/4a^2+ab+b^2+1/2ac+bc+3/4c^2=3/4$, which has 8 integer solutions. We know $\sum_{v\in {R}}{contr_{v}(P_3)}=2\chi(S)+2(O.P_3)-<P_3,P_3>$ and $max(\sum_{v\in {R}}{contr_{v}(P_3))=5/4},$ so $2+2(O.P_3)-1\leq 5/4$, and $(O.P_3)\leq 1/8$. We have $(O.P_3)=0$. In fact, we have 8 integral sections  as follows:

\begin{itemize}
  \item (0,0,$\pm1$), ($\pm1$,0,0), (0,1,-1), (0,-1,1), (1,-1,0), (-1,1,0)
\end{itemize}

This computation method can be extended directly to other case. Here all cases can be treated by exactly the same arguments as above.
% Order to avoid the proof method being too repetitive, so next, for the same proof method, we will briefly discuss the following cases, and for special situation, we will make special analysis.

{\bf Proof of  Case $\uppercase\expandafter{\romannumeral 2}$}

Coming to the case of second category above, we only touch upon a  relevant case No. 24.  $T=A^5_1$, where $v(T)=10$ and Mordell-Weil group $E(K) \cong A^{\vee}_3 \bigoplus \mathbb{Z}_2$.  There are 5  reducible fibres of type $\uppercase\expandafter{\romannumeral 3}$.

 By the explicit formula of the height pairing, we have $P \in E(K)$ as follows:
$$<P, P>= 2 + 2(P.O)-\left\{
\begin{array}{ll}
0, & \hbox{$P\in E(K)^0$, $P$ intersects $\Theta_{v_k,0}~(k=1,2,3,4,5)$;} \\
5/2, & \hbox{$P\in E(K)-E(K)^0$, $P$ intersects $\Theta_{v_k,1}~(k=1,2,3,4,5)$;} \\
2, & \hbox{$P\in E(K)-E(K)^0$, $P$ intersects $\Theta_{v_k,1}~(k=1,2,3,4)$;} \\
3/2, & \hbox{$P\in E(K)-E(K)^0$, $P$ intersects $\Theta_{v_k,1}~(k=1,2,3)$;}\\
1, & \hbox{$P\in E(K)-E(K)^0$, $P$ intersects $\Theta_{v_k,1}~(k=1,2)$;} \\
1/2, & \hbox{$P\in E(K)-E(K)^0$, $P$ intersects $\Theta_{v_1,1}$.}
\end{array}
\right.
$$

$P\in E(K)-E(K)^0$ and $\sum_{v\in {R}}{contr_{v}(P)}=5/2$, so $<P,P>=-1/6$. However, $<P, P> = (\varphi(P).\varphi(P))$ this defines a $\mathbb{Q}$-valued symmetric bilinear pairing on $E(K)$ which induces the structure of a positive-definite lattice on $E(K)/E(K)_{tors}$. Hence, there is no such $P$. If assume that $P_i \in {\mathcal{P}_i}$, thus we get $<P_i, P_i>=2,~0,~1/2,~1,~3/2~(i=1,~2,~3,~4,~5)$.

First, according to different the height pairing value of sections, we compute  their  $m(P_i)$. Here we can use the same method from the previous case.
Now we have ${D_p}^2=-2$, then $((P)-(O)-(F)+\sum_{{v_k}, i}{n_{{v_k},i}}\Theta_{{v_k},i})^2=-2 ~(i=\{1\},~ k=\{1, \cdots, 5\})$ with $(P.O)=0, ~F^2=0, ~(P.P)=(O,O)=-\chi(S)=-1, ~(\Theta_{{v_k},1})^2=-2$. Since it has 5 singular fibres $\uppercase\expandafter{\romannumeral 3}$, then $(\Theta_{{v_k},1}.\Theta_{{v_{k'}},1})=0, k\neq k' $.  Moreover, according to  different intersection ways of $P$ and reducible fibre's components, intersection number of $P$ and $\Theta_{{v_k},1}$ is 1, so we obtain 5 equations as follows.
\begin{equation}
\begin{aligned}
(1): &\sum\limits_{i=1}^5 (n_{{v_i},1})^2=0;~~(2): \sum\limits_{i=1}^5 (n_{{v_i},1})^2-n_{{v_1},1}-n_{{v_2},1}-n_{{v_3},1}-n_{{v_4},1}=0;\\
(3): &\sum\limits_{i=1}^5 (n_{{v_i},1})^2-n_{{v_1},1}-n_{{v_2},1}-n_{{v_3},1}=0;~~(4): \sum\limits_{i=1}^5 (n_{{v_i},1})^2-n_{{v_1},1}-n_{{v_2},1}=0;\\
(5): &\sum\limits_{i=1}^5 (n_{{v_i},1})^2-n_{{v_1},1}=0.
\end{aligned}
\end{equation}
After we compute  the non-negative integer solutions of these five equations, then it is clear that $<P_i, P_i>=2,~0,~1/2,~1,~3/2$, hence, $m(P_i)=1,~16, ~8, ~4, ~2 ~(i=1,~2,~3,~4,~5)$.

Second, we compute the number of their integral sections. In the same way of the case No.16, when $<P_1, P_1>=2$, we have ${\mathrm \#} {\mathcal{P}_1}=6$.
When $<P_2, P_2>=0$, because $E(K)$ has a symmetric
bilinear pairing which is positive-definite modulo torsion, we derive the equivalence $<P_2, P_2>=0$ iff $\varphi (P_2)=0$ iff $P_2\in E(K)_{tor}$, hence $P_2 \in \mathbb{Z}_2$. From \cite{SS}, we know that for any torsion element $P$ of order not divisible by char($k$), the sections $(P)$ and $(O)$ are disjoint. Here $char(k)=0$,  it is clear that  it has only 1 non-trivial integral section (not contain section $(O)$).

 %In the same way of above,  since intersection number of $P$ and $\Theta_{{v_k},1}$ is 1 and $P$ intersect components $\Theta_{v_i,1} (i=1,2,3,4)$ . We get equation $(\ref{equation1})-\sum_{k}({n_{{v_k},1}})^2=0$, this equation
%can be computed number of solutions, which have 16 integral  solutions, it means  $m(P_2)=16$. In the following, we can use the same method to get $m(P_3)=8, m(P_4)=4$ and $m(P_5)=2$

When $<P_3, P_3>=1/2, ~P_3\in E(K)-E(K)^0$, we get ${\mathrm \#} {\mathcal{P}_3}=6$. Suppose a torsion element $\alpha \in {E(K)_{tor}}$, then $<P_3, P_3>=<P_3+\alpha, P_3+\alpha>=1/2$ by Lemma \ref{lemma2}, here $\alpha=\bar{1}\in \mathbb{Z}_2$. According to the height pairing formula. It is easy to check that following result: if $P_3$ is an integral section, then $P_3+\alpha$ is also an integral section. Finally, we can obtain 12 integral sections. When $<P_4, P_4>=1$, it has 24 integral sections in the same way.

When $<P_5, P_5>=3/2, ~P_5\in E(K)-E(K)^0$. Here, we find that it has Situation $(S2)$. If assume that $P_5$ is an integral section, we get ${\mathrm \#} P_5=8$. However, while we are checking whether it is integral section,  observe that $(O.P_5)\leq 1$ by the height pairing formula. It means that the number of integral sections less than or equal to 16 by Lemma \ref{lemma2}, hence, we can assume that number of integral section is ${\mathrm \#} {\mathcal{P}_5} :=x$. We have
$240-v(T)=1\cdot16+6\cdot1+24\cdot4+12\cdot8+2x$,  this can easily be obtained from Theorem \ref{theorem5}, so $x=8$.

{\bf Proof of case $\uppercase\expandafter{\romannumeral 3}$}

First, we list an instance and turn to case No. 51 in case $\uppercase\expandafter{\romannumeral 3}$.  $T=A_5\bigoplus A_2$, $v(T)=36$ and Mordell-Weil group $E(K) \cong A^{\vee}_1 \bigoplus \mathbb{Z}_3$. The two reducible fibres are $I_6$ and $I_3$.
By the explicit formula of the height pairing, we have for $P \in E(K)$, then
$$ <P, P>= 2 + 2(P.O)-\left\{
\begin{array}{ll}
0, & \hbox{$P\in E(K)^0$, $P$ intersects $\Theta_{v_1,0}$ and $\Theta_{v_2,0}$;} \\
2, & \hbox{$P\in E(K)-E(K)^0$, $P$ intersects $\Theta_{{v_1},2}$ and $\Theta_{{v_2},1}$;} \\
5/6, & \hbox{$P\in E(K)-E(K)^0$, $P$ intersects $\Theta_{v_1,1}$ and $\Theta_{{v_2},0}$;}\\
4/3, & \hbox{$P\in E(K)-E(K)^0$, $P$ intersects $\Theta_{v_1,2}$ and $\Theta_{{v_2},0}$;}\\
3/2, & \hbox{$P\in E(K)-E(K)^0$, $P$ intersects $\Theta_{v_1,3}$ and $\Theta_{{v_2},0}$;} \\
3/2, & \hbox{$P\in E(K)-E(K)^0$, $P$ intersects $\Theta_{v_1,1}$ and $\Theta_{{v_2},1}$;} \\
2/3, & \hbox{$P\in E(K)-E(K)^0$, $P$ intersects $\Theta_{v_1,0}$ and $\Theta_{{v_2},1}$;} \\
13/6, & \hbox{$P\in E(K)-E(K)^0$, $P$ intersects $\Theta_{v_1,3}$ and $\Theta_{{v_2},1}$.}
\end{array}
\right.
$$

There is no such $P$ such that $\sum_{v\in {R}}{contr_{v}(P)}=13/6$.
If we assume that $P_i ~(i=1,~2,~3,~4,~5,~6,~7)$ are integral sections, then we get $<P_i, P_i>=2, ~0, ~7/6, ~2/3, ~1/2, ~1/2, ~4/3$.
When $<P_i, P_i>=7/6, ~2/3, ~4/3 ~(i=3, ~4, ~7)$, it has  $P\in E(K)-E(K)^0$. Assume that $e_1$ is a basis of $E(K)$ such that  $P=a~e_1$ by the Gram matrix of $A^{\vee}_1$, it is easy to get $1/2a^2=7/6, ~2/3, ~4/3$, which have no integer solution.

First, we have $m(P_i)=1, ~45 ~~(i=1, ~2)$ and ${\mathrm \#} \mathcal{P}_1=2, ~{\mathrm \#} \mathcal{P}_2=2$ by an easy computation.
When $<P_i, P_i>=1/2~(i=5,6)$,  we get 6 integral sections. But in this situation, we know that it just has Situation $(S2)$. The first is that $P_5$ intersects components $\Theta_{v_1,3}$ and $\Theta_{{v_2},0}$, and corresponding to $m(P_5)=20$, we assume that the number of integral sections
 is ${\mathrm \#} \mathcal{P}_5 :=x$. The second situation is that $P$ intersects components $\Theta_{v_1,1}$ and $\Theta_{{v_2},1}$, and corresponding to $m(P_6)=18$, we assume that the number of integral section is ${\mathrm \#} \mathcal{P}_6 :=y$, $x+y=6$. In fact, we can then apply the same method such that we have $240-v(T)=1\cdot2+45\cdot2+20x+18y$, so $x=2, ~y=4$.

{\bf Proof of case $\uppercase\expandafter{\romannumeral 4}$}

This situation is slightly more complicated than the previous three cases.  Finally, as an explicit illustration, we will go into details only for an instructive case: No.53.

\begin{figure}[ht]
\begin{center}
\scalebox{0.8}{\includegraphics{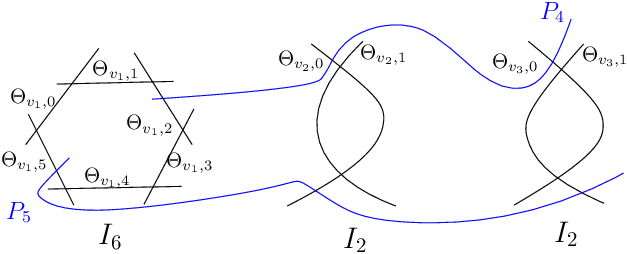}}
\end{center} \caption{}\label{No53}
\end{figure}

We consider the case No.53 where $T=A_5\bigoplus A^2_1$, $v(T)=34$ and  $E(K) \cong <1/6> \bigoplus \mathbb{Z}_2$.  The three reducible fibres are $I_6$, $I_2$ and $I_2$. As discussed above, we get $$\sum_{v\in {R}}{contr_{v}(P)}=0, ~2, ~5/6, ~4/3, ~4/3, ~3/2, ~1/2, ~1, ~11/6, ~11/6, ~3/7, ~5/2.$$ When $\sum_{v\in {R}}{contr_{v}(P)}=3/7, ~5/2$, there is no such $P$. It is clear that $$<P_i, P_i>=2, ~0, ~7/6, ~2/3, ~2/3, ~1/2, ~3/2, ~1, ~1/6, ~1/6~ (i=1, ~2, ~3, ~4, ~5, ~6, ~7, ~8, ~9, ~10).$$ When $<P_i, P_i>=2, ~7/6, ~1/2, ~1~(i=1, ~3, ~6, ~8)$, which have no integer solution. We assume that $P_i ~(i=2, ~4, ~5, ~7, ~9, ~10)$ are integral sections.

When $<P_2, P_2>=0$, we can proceed similarly for the computation, it gives rise to ${\mathrm \#} {\mathcal{P}_1}=1$ and $m(P_1)=40$.

When $<P_7, P_7>=3/2$  in this situation, we know that it just has Situation ($S1$). It has  $m(P_7)=2$. But it gives rise to $(O.P_7)\leq 1$ by checking whether it is an integral section or not. Hence, we can assume that the number of such integral sections is ${\mathrm \#} \mathcal{P}_7:=x$. As we remarked in  Lemma \ref{lemma2}, it proves $x\leq4$.

When $<P_i, P_i>=2/3 ~(i=4, ~5)$, it has 4 integral sections in all.
But in this situation, we know that it just has Situation ($S2$).  The first  is that $P_4$ should pass through the  components $\Theta_{v_1,2}$, $\Theta_{{v_2},0}$ and $\Theta_{{v_3},0}$ of every reducible fibre~(as depicted in Figure \ref{No53} )~and corresponding to $m(P_4)=15$. Assume that the number of such integral sections is ${\mathrm \#} \mathcal{P}_4:=y$. The second situation is that $P_5$ intersects components $\Theta_{v_1,5}$, ~$\Theta_{{v_2},0}$ and $\Theta_{{v_3},1}$, and corresponding to $m(P_5)=12$, we assume that the number of such integral sections is ${\mathrm \#} \mathcal{P}_5:=z$. In fact $y+z=4$.

When $<P_i, P_i>=1/6~(i=9, ~10)$, it has 4 integral sections in all. And
it just has Situation ($S2$) as above.  The first  is that $P_9$ intersects components $\Theta_{v_1,4}$, $\Theta_{{v_2},1}$ and $\Theta_{{v_3},0}$, and corresponding to $m(P_9)=30$. We assume that the number of such integral sections
is ${\mathrm \#} \mathcal{P}_9:=w$. The second situation is that $P_{10}$ intersects components $\Theta_{v_1,1}$, $\Theta_{{v_2},1}$ and $\Theta_{{v_3},1}$, and corresponding to $m(P_{10})=24$, we assume that the number of such integral sections is ${\mathrm \#} \mathcal{P}_{10}:=t$. In fact $w+t=4$.

In a word, according to relevant computation value of $P_i$, we have $240-v(T)=40\cdot1+2x+15y+12z+30w+24t$.
The non-negative integer solutions $(x,~y,~z,~w,~t)$ satisfying the conditions are (2,~4,~0,~1,~3),~(2,~2,~2,~2,~2),~(2,~0,~4,~3,~1).  We know that
$O \neq P_1 \in E(K)_{tor}$, and $P_1$ intersects components $\Theta_{v_1,3}$, $\Theta_{{v_2},0}$ and $\Theta_{{v_3},1}$. It is easy to check  $P_4+P_1 \in \mathcal{P}_5$, $P_5+P_1 \in \mathcal{P}_4$, $P_9+P_1 \in \mathcal{P}_{10}$ and $P_{10}+P_1 \in \mathcal{P}_9$.  It implies that
$y=z$ and $w=t$, hence the unique solution is (2,~2,~2,~2,~2).
The other cases in this situation are similar. This completes the proof of Theorem \ref{thm1}.

\section{Proof of   Theorem 2}\label{proof1}
In this section, we will prove Theorem \ref{thm2}.
First, we recall a lemma in \cite{S2}. It helps us  estimate the number of  solutions of the Weierstrass  equation as follows:
$$ y^2=x^3+a_2(t)x^2+a_4(t)x+a_6(t), ~\mathrm{deg} a_i(t))\leq i ~~~~~~ (*)$$

\begin{lemma}
Let $P =(x, y) \in E(K), P \neq O$. Then the section
$(P)$ is an integral section, if and only if $x$ and $y$ are
polynomials for $t$ of degree $\leq 2$ or $\leq 3$, i.e., of the form:
$$ x = gt^2 + at + b,~~~~y = ht^3 + ct^2 + dt + e$$
with $a, b,\ldots, g, h \in k.$
\end{lemma}

Now, we start our proof.
\begin{proof}

(1)~~If $\Delta(t)$ has only one solution which is not 12-fold root. Then the Weierstrass equation $(*)$  is minimal. And it must correspond to an elliptic fibration over $\mathbb{P}^1$.  This fibration has exactly two singular fibers by \cite{Bea}. Then by \cite{Hiz} or  \cite{Ngu}, its $T$ must be one of $E_8, E_7\oplus A_1, E_6\oplus A_2, D_4\oplus D_4$. In these cases,  every section is an integral section and corresponds to a  polynomial solution (see \cite{TI}).  From the table in Theorem \ref{thm1},  we have $0\leq n\leq3$.

(2)~~If $\Delta(t)$ has only two solutions or  three different solutions with deg $\Delta(t) = 12$. Then $(*)$ corresponds to an elliptic fibration over $\mathbb{P}^1$ with two or three singular fibers. From the classifications in \cite{Hiz} or \cite{Ngu}, its $T$ must be one of the followings:
$$E_8, E_7\oplus A_1, E_6\oplus A_2, D_4\oplus D_4, D_5\oplus A_3, D_6\oplus {A_1}^2, D_8, D_6\oplus A_1,$$ $$D_5\oplus {A_1}^2, D_5\oplus A_2, E_7, E_6\oplus A_1, D_7, {A_2}^3, D_4\oplus {A_1}^2, D_5\oplus A_1, D_6, D_4\oplus A_2, E_6.$$  From the table in Theorem \ref{thm1},  we  have $ n\leq N_0$, and $N_0$ is at least 26. (In fact, we know that $N_0\leq 144\times10^{14.2}$ from \cite{HS}.)
\begin{itemize}
     \item (3)-(i)~~If $\Delta(t)$ has more than  three solutions, or three different solutions with deg $\Delta(t) < 12$. Then it corresponds to an elliptic fibration over $\mathbb{P}^1$ with four singular fibers.
We compute Euler-Poincare characteristic $\chi_{top}{(F)}$ for every elliptic fiber in the following:

\begin{itemize}
  \item $Fiber~type$ :   $I_n~(n\geq0)$~~~${I_n}^*~(n\geq0)$~~~{\uppercase\expandafter{\romannumeral 2}}$~~~${\uppercase\expandafter{\romannumeral 3}}$~~~${\uppercase\expandafter{\romannumeral 4}}~~~${\uppercase\expandafter{\romannumeral 2}}^{*}$~~~${\uppercase\expandafter{\romannumeral 3}}^{*}$ ~~~${\uppercase\expandafter{\romannumeral 4}}^{*}$
  \item $\chi_{top}{(F)}$ : ~~~~~~~~~~~~$n$~~~~~~~~~$6+n$~~~~~~~$2$~~~~~$3$~~~~~$4$~~~~~$10$~~~~~~~$9$~~~~~~~$8$
\end{itemize}
We have the formula $\chi_{top}(S)=\sum_{i=1}^{n} \chi_{top}(F)$. Because $S$ is a rational surface, $\chi_{top}{(S)}=12$.
If an elliptic fibration with $II^{*}$ and other three singular fibers,  then $\chi_{top}{(S)}\geq 10+1+1+1$.  This is a contradiction.  It means that there are no fibration with ${\uppercase\expandafter{\romannumeral 2}}^{*}$ in this case. From the table in Theorem \ref{thm1},  and we have $n\geq 1$.
     \item (3)-(ii)~~If $3 \leq  \mathrm{deg} \Delta(t) \leq9$, it means  the valuation of $\Delta(t)$  must be not $10$ at every point in $\mathbb{P}^1$.  From Tate's algorithm, there is no fibration with ${\uppercase\expandafter{\romannumeral 2}}^{*}$ in this case.
   \end{itemize}
\end{proof}

Now, in the following we give some examples. These examples show that the bounds are attainable  in Theorem 2.
\\{\bf Example~~9}

\begin{enumerate}
  \item The Weiestrass equation $y^2=x^3+t^3$ corresponds to an ellptic fibration over $\mathbb{P}^1$ with ${I_0}^{*} (t=\infty) $ and ${I_0}^{*}$ (t=0). We compute  $\Delta(t)=-27t^6$. It is easy to see  $T=D_4 \bigoplus D_4 $. From the table in Theorem \ref{thm1},  we have $n=3$.
  \item The Weiestrass equation $y^2=x^3-3tx^2+3t^2x+t^2+t^3+t^4$ corresponds to an ellptic fibration over $\mathbb{P}^1$ with $IV$, $IV$  and $IV$, at ~$t=0, ~t=-1, ~t=\infty$. It is easy to see  $T= A_2^3$. We  compute $\Delta(t)=-27t^4(1+t)^4$. From the table,  we have $n \geq 26$.
  \item The Weiestrass equation $y^2=x^3+t^2x^2+(t^3+t^4)x+t^6$ corresponds to an elliptic fibration over $\mathbb{P}^1$ with  $III^{*}$ (at t=0), and $I_1,I_1, I_1$. It is easy to see $T= E_7$. We compute $\Delta(t)=-16t^{12}+8t^{11}-11t^{10}-4t^9$. We can obtain $n=1$ by the table in Theorem \ref{thm1}.
\end{enumerate}

\begin{table}[ht]
    \centering
    \footnotesize
    \renewcommand{\arraystretch}{1.0}
    \setlength{\tabcolsep}{10pt}
    \caption{Main table}
    \label{tab1:my_label}
    \scalebox{1.1}{
\begin{tabular}{ccccccc}
\hline
No. &  $r$  & T                                     & $v(T)$ & $<P_i,P_i>$ & $m(P_i)$ & ${\mathrm \#} {\mathcal{P}_i}$ \\ \hline
1   & 8    & \{0\}                                 & 0    & 2                            & 1    & 240 \\ \hline
2   & 7    & $A_1$                                 & 2    & 2                            & 1    & 126 \\ \hline
    &      &                                       &      & 3/2                          & 2    & 56  \\ \hline
3   & 6    & $A_2$                                 & 6    & 2                            & 1    & 72  \\ \hline
    &      &                                       &      & 4/3                          & 3    & 54  \\ \hline
4   & 6    & ${A_1}^2$                             & 4    & 2                            & 1    & 60  \\ \hline
    &      &                                       &      & 3/2                          & 2    & 64  \\ \hline
    &      &                                       &      & 1                            & 4    & 12  \\ \hline
5   & 5    & $A_3$                                 & 12   & 2                            & 1    & 40  \\ \hline
    &      &                                       &      & 5/4                          & 4    & 32  \\ \hline
    &      &                                       &      & 1                            & 6    & 10  \\ \hline
6   & 5    & $A_2 \bigoplus A_1$                   & 8    & 2                            & 1    & 30  \\ \hline
    &      &                                       &      & 3/2                          & 2    & 20  \\ \hline
    &      &                                       &      & 4/3                          & 3    & 30  \\ \hline
    &      &                                       &      & 5/6                          & 6    & 12  \\ \hline
7   & 5    & ${A_1}^3$                             & 6    & 2                            & 1    & 26  \\ \hline
    &      &                                       &      & 3/2                          & 2    & 48  \\ \hline
    &      &                                       &      & 1                            & 4    & 24  \\ \hline
    &      &                                       &      & 1/2                          & 8    & 2   \\ \hline
8   & 4    & $A_4$                                 & 20   & 2                            & 1    & 20  \\ \hline
    &      &                                       &      & 6/5                          & 5    & 20  \\ \hline
    &      &                                       &      & 4/5                          & 10   & 10  \\ \hline
9   & 4    & $D_4$                                 & 24   & 2                            & 1    & 24  \\ \hline
    &      &                                       &      & 1                            & 8    & 24  \\ \hline
10  & 4    & $A_3 \bigoplus A_1$                   & 14   & 2                            & 1    & 14  \\ \hline
    &      &                                       &      & 3/2                          & 2    & 12  \\ \hline
    &      &                                       &      & 5/4                          & 4    & 16  \\ \hline
    &      &                                       &      & 1                            & 6    & 6   \\ \hline
    &      &                                       &      & 3/4                          & 8    & 8   \\ \hline
    &      &                                       &      & 1/2                          & 12   & 2   \\ \hline
\end{tabular}
}
\end{table}

\begin{table}[ht]
    \centering
    \footnotesize
    \renewcommand{\arraystretch}{1.0}
    \setlength{\tabcolsep}{10pt}
    \caption{Main table}
    \label{tab2:my_label}
    \scalebox{1.1}{
\begin{tabular}{ccccccc}
\hline
No. & $r$  & T                                     & $v(T)$ & $<P_i,P_i>$ & $m(P_i)$ & ${\mathrm \#} {\mathcal{P}_i}$ \\ \hline
11  & 4    & ${A_2}^2$                             & 12   & 2                            & 1    & 12  \\ \hline
    &      &                                       &      & 4/3                          & 3    & 36  \\ \hline
    &      &                                       &      & 2/3                          & 9    & 12  \\ \hline
12  & 4    & $A_2 \bigoplus {A_1}^2$               & 10   & 2                            & 1    & 12  \\ \hline
    &      &                                       &      & 3/2                          & 2    & 16  \\ \hline
    &      &                                       &      & 4/3                          & 3    & 14  \\ \hline
    &      &                                       &      & 1                            & 4    & 6   \\ \hline
    &      &                                       &      & 5/6                          & 6    & 16  \\ \hline
    &      &                                       &      & 1/3                          & 12   & 2   \\  \hline
13  & 4    & ${A_1}^4$                             & 8    & 0                            & 16   & 1   \\ \hline
    &      &                                       &      & 2                            & 1    & 24  \\ \hline
    &      &                                       &      & 1                            & 4    & 48  \\ \hline
14  & 4    & ${A_1}^4$                             & 8    & 2                            & 1    & 8   \\ \hline
    &      &                                       &      & 3/2                          & 2    & 32  \\ \hline
    &      &                                       &      & 1                            & 4    & 24  \\ \hline
    &      &                                       &      & 1/2                          & 8    & 8   \\ \hline
15  & 3    & $A_5$                                 & 30   & 2                            & 1    & 8   \\ \hline
    &      &                                       &      & 7/6                          & 6    & 12  \\ \hline
    &      &                                       &      & 2/3                          & 15   & 6   \\ \hline
    &      &                                       &      & 1/2                          & 20   & 2   \\ \hline
16  & 3    & $D_5$                                 & 40   & 2                            & 1    & 12  \\ \hline
    &      &                                       &      & 1                            & 10   & 6   \\ \hline
    &      &                                       &      & 3/4                          & 16   & 8   \\ \hline
17  & 3    & $A_4 \bigoplus A_1$                   & 22   & 2                            & 1    & 6   \\ \hline
    &      &                                       &      & 3/2                          & 2    & 6   \\ \hline
    &      &                                       &      & 6/5                          & 5    & 8   \\ \hline
    &      &                                       &      & 4/5                          & 10   & 6   \\ \hline
    &      &                                       &      & 7/10                         & 10   & 6   \\ \hline
    &      &                                       &      & 3/10                         & 20   & 2   \\ \hline
\end{tabular}
}
\end{table}

\begin{table}[ht]
    \centering
    \footnotesize
    \renewcommand{\arraystretch}{1.0}
    \setlength{\tabcolsep}{10pt}
    \caption{Main table}
    \label{tab3:my_label}
    \scalebox{1.1}{
\begin{tabular}{ccccccc}
\hline
No. & $r$  & T                                     & $v(T)$ & $<P_i,P_i>$ & $m(P_i)$ & ${\mathrm \#} {\mathcal{P}_i}$ \\ \hline
18  & 3    & $D_4 \bigoplus A_1$                   & 26   & 2                            & 1    & 6   \\ \hline
    &      &                                       &      & 3/2                          & 2    & 8   \\ \hline
    &      &                                       &      & 1                            & 8    & 12  \\ \hline
    &      &                                       &      & 1/2                          & 16   & 6   \\ \hline
19  & 3    & $A_3 \bigoplus A_2$                   & 18   & 2                            & 1    & 4   \\ \hline
    &      &                                       &      & 4/3                          & 3    & 10  \\ \hline
    &      &                                       &      & 5/4                          & 4    & 8   \\ \hline
    &      &                                       &      & 1                            & 6    & 4   \\ \hline
    &      &                                       &      & 7/12                         & 12   & 8   \\ \hline
    &      &                                       &      & 1/3                          & 18   & 2   \\ \hline
20  & 3    & ${A_2}^2 \bigoplus A_1$               & 14   & 2                            & 1    & 6   \\ \hline
    &      &                                       &      & 3/2                          & 2    & 2   \\ \hline
    &      &                                       &      & 4/3                          & 3    & 12  \\ \hline
    &      &                                       &      & 5/6                          & 6    & 12  \\ \hline
    &      &                                       &      & 2/3                          & 9    & 8   \\ \hline
    &      &                                       &      & 1/6                          & 18   & 2   \\ \hline
21  & 3    & $A_3 \bigoplus {A_1}^2$               & 16   & 0                            & 24   & 1   \\ \hline
    &      &                                       &      & 2                            & 1    & 12  \\ \hline
    &      &                                       &      & 1                            & 6    & 6   \\  \hline
    &      &                                       &      & 1                            & 4    & 6   \\  \hline
    &      &                                       &      & 3/4                          & 8    & 16  \\  \hline
22  & 3    & $A_3 \bigoplus {A_1}^2$               & 16   & 2                            & 1    & 4   \\  \hline
    &      &                                       &      & 3/2                          & 2    & 8   \\ \hline
    &      &                                       &      & 5/4                          & 4    & 8   \\ \hline
    &      &                                       &      & 1                            & 6    & 2   \\ \hline
    &      &                                       &      & 1                             & 4    & 4   \\ \hline
    &      &                                       &      & 3/4                          & 8    & 8   \\ \hline
    &      &                                       &      & 1/2                          & 12   & 4   \\ \hline
    &      &                                       &      & 1/4                          & 16   & 2   \\ \hline
\end{tabular}
}
\end{table}

\begin{table}[ht]
    \centering
    \footnotesize
    \renewcommand{\arraystretch}{1.0}
    \setlength{\tabcolsep}{10pt}
    \caption{Main table}
    \label{tab4:my_label}
    \scalebox{1.1}{
\begin{tabular}{ccccccc}
\hline
No. & $r$  & T                                     & $v(T)$ & $<P_i,P_i>$ & $m(P_i)$ & ${\mathrm \#} {\mathcal{P}_i}$ \\ \hline
23  & 3    & $A_2 \bigoplus {A_1}^3$               & 12   & 2                            & 1    & 2   \\ \hline
    &      &                                       &      & 3/2                          & 2    & 12  \\ \hline
    &      &                                       &      & 4/3                          & 3    & 6   \\ \hline
    &      &                                       &      & 1                            & 4    & 6   \\ \hline
    &      &                                       &      & 5/6                          & 6    & 12  \\ \hline
    &      &                                       &      & 1/2                          & 8    & 2   \\ \hline
    &      &                                       &      & 1/3                          & 12   & 6   \\  \hline
24  & 3    & ${A_1}^5$                             & 10   & 0                            & 16   & 1   \\ \hline
    &      &                                       &      & 2                            & 1    & 6   \\ \hline
    &      &                                       &      & 3/2                          & 2    & 8   \\ \hline
    &      &                                       &      & 1                            & 4    & 24  \\ \hline
    &      &                                       &      & 1/2                          & 8    & 12  \\ \hline
25  & 2    & $A_6$                                 & 42   & 2                            & 1    & 2   \\ \hline
    &      &                                       &      & 8/7                          & 7    & 6   \\ \hline
    &      &                                       &      & 4/7                          & 21   & 4   \\ \hline
    &      &                                       &      & 2/7                          & 35   & 2   \\ \hline
26  & 2    & $D_6$                                 & 60   & 2                            & 1    & 4   \\ \hline
    &      &                                       &      & 1                            & 12   & 4   \\ \hline
    &      &                                       &      & 1/2                          & 32   & 4   \\ \hline
27  & 2    & $E_6$                                 & 72   & 2                            & 1    & 6   \\ \hline
    &      &                                       &      & 2/3                          & 27   & 6   \\ \hline
28  & 2    & $A_5 \bigoplus A_1$                   & 32   & 0                            & 40   & 1   \\ \hline
    &      &                                       &      & 2                            & 1    & 6   \\ \hline
    &      &                                       &      & 2/3                          & 15   & 6   \\ \hline
    &      &                                       &      & 2/3                          & 12   & 6   \\ \hline
\end{tabular}
}
\end{table}

\begin{table}[ht]
    \centering
    \footnotesize
    \renewcommand{\arraystretch}{1.0}
    \setlength{\tabcolsep}{10pt}
    \caption{Main table}
    \label{tab5:my_label}
    \scalebox{1.1}{
\begin{tabular}{ccccccc}
\hline
No. & $r$  & T                                     & $v(T)$ & $<P_i,P_i>$ & $m(P_i)$ & ${\mathrm \#} {\mathcal{P}_i}$ \\ \hline
29  & 2    & $A_5 \bigoplus A_1$                   & 32   & 2                            & 1    & 2   \\ \hline
    &      &                                       &      & 3/2                          & 2    & 2   \\ \hline
    &      &                                       &      & 7/6                          & 6    & 4   \\ \hline
    &      &                                       &      & 2/3                          & 15   & 2   \\ \hline
    &      &                                       &      & 2/3                          & 12   & 4   \\ \hline
    &      &                                       &      & 1/2                          & 20   & 2   \\ \hline
    &      &                                       &      & 1/6                          & 30   & 2   \\ \hline
30  & 2    & $D_5 \bigoplus A_1$                   & 42   & 2                            & 1    & 2   \\ \hline
    &      &                                       &      & 3/2                          & 2    & 4   \\ \hline
    &      &                                       &      & 1                            & 10   & 2   \\ \hline
    &      &                                       &      & 3/4                          & 16   & 4   \\ \hline
    &      &                                       &      & 1/2                          & 20   & 2   \\ \hline
    &      &                                       &      & 1/4                          & 32   & 2   \\ \hline
31  & 2    & $A_4 \bigoplus A_2$                   & 26   & 2                            & 1    & 2   \\ \hline
    &      &                                       &      & 6/5                          & 5    & 2   \\ \hline
    &      &                                       &      & 4/3                          & 3    & 4   \\ \hline
    &      &                                       &      & 4/5                          & 10   & 4   \\ \hline
    &      &                                       &      & 8/15                         & 15   & 6   \\ \hline
    &      &                                       &      & 2/15                         & 30   & 2   \\ \hline
32  & 2    & $D_4 \bigoplus A_2$                   & 30   & 4/3                          & 3    & 6   \\ \hline
    &      &                                       &      & 1                            & 8    & 6   \\ \hline
    &      &                                       &      & 1/3                          & 24   & 6   \\ \hline
33  & 2    & $A_4 \bigoplus {A_1}^2$               & 24   & 3/2                          & 2    & 4   \\ \hline
    &      &                                       &      & 6/5                          & 5    & 4   \\ \hline
    &      &                                       &      & 1                            & 4    & 2   \\ \hline
    &      &                                       &      & 4/5                          & 10   & 2   \\ \hline
    &      &                                       &      & 7/10                         & 10   & 4   \\ \hline
    &      &                                       &      & 3/10                         & 20   & 4   \\ \hline
    &      &                                       &      & 1/5                          & 20   & 2   \\ \hline
\end{tabular}
}
\end{table}

\begin{table}[ht]
    \centering
    \footnotesize
    \renewcommand{\arraystretch}{1.0}
    \setlength{\tabcolsep}{10pt}
    \caption{Main table}
    \label{tab6:my_label}
    \scalebox{1.1}{
\begin{tabular}{ccccccc}
\hline
No. & $r$  & T                                     & $v(T)$ & $<P_i,P_i>$ & $m(P_i)$ & ${\mathrm \#} {\mathcal{P}_i}$ \\ \hline
34  & 2    & $D_4 \bigoplus {A_1}^2$               & 28   & 0                            & 32   & 1   \\ \hline
    &      &                                       &      & 2                            & 1    & 4   \\ \hline
    &      &                                       &      & 1                            & 8    & 4   \\ \hline
    &      &                                       &      & 1                            & 4    & 4   \\ \hline
    &      &                                       &      & 1/2                          & 16   & 8   \\ \hline
35  & 2    & ${A_3}^2$                             & 24   & 0                            & 36   & 1   \\ \hline
    &      &                                       &      & 2                            & 1    & 4   \\ \hline
    &      &                                       &      & 1                            & 6    & 8   \\ \hline
    &      &                                       &      & 1/2                          & 16   & 8  \\ \hline
36  & 2    & ${A_3}^2$                             & 24   & 5/4                          & 4    & 8   \\ \hline
    &      &                                       &      & 1                            & 6    & 4   \\ \hline
    &      &                                       &      & 1/2                          & 16   & 4   \\ \hline
    &      &                                       &      & 1/4                          & 24   & 4   \\ \hline
37  & 2    & $A_3 \bigoplus A_2 \bigoplus A_1$     & 20   & 2                            & 1    & 2   \\ \hline
    &      &                                       &      & 4/3                          & 3    & 2   \\ \hline
    &      &                                       &      & 5/4                          & 4    & 4   \\ \hline
    &      &                                       &      & 5/6                          & 6    & 4   \\ \hline
    &      &                                       &      & 3/4                          & 8    & 2   \\ \hline
    &      &                                       &      & 7/12                         & 12   & 4   \\ \hline
    &      &                                       &      & 1/2                          & 12   & 2   \\ \hline
    &      &                                       &      & 1/3                          & 18   & 2   \\ \hline
    &      &                                       &      & 1/12                         & 24   & 2   \\ \hline
38  & 2    & $A_3 \bigoplus {A_1}^3$               & 18   & 0                            & 24   & 1   \\ \hline
    &      &                                       &      & 2                            & 1    & 2   \\ \hline
    &      &                                       &      & 3/2                          & 2    & 4   \\ \hline
    &      &                                       &      & 1                            & 6    & 2   \\ \hline
    &      &                                       &      & 1                            & 4    & 2   \\ \hline
    &      &                                       &      & 3/4                          & 8    & 8   \\ \hline
    &      &                                       &      & 1/2                          & 12   & 2   \\ \hline
    &      &                                       &      & 1/2                          & 8    & 2   \\ \hline
    &      &                                       &      & 1/4                          & 16   & 4   \\ \hline
\end{tabular}
}
\end{table}

\begin{table}[ht]
    \centering
    \footnotesize
    \renewcommand{\arraystretch}{1.0}
    \setlength{\tabcolsep}{10pt}
    \caption{Main table}
    \label{tab7:my_label}
    \scalebox{1.1}{
\begin{tabular}{ccccccc}
\hline
No. & $r$  & T                                     & $v(T)$ & $<P_i,P_i>$ & $m(P_i)$ & ${\mathrm \#} {\mathcal{P}_i}$ \\ \hline
39  & 2    & ${A_2}^3$                             & 18   & 0                            & 27   & 2   \\ \hline
    &      &                                       &      & 2                            & 1    & 6   \\ \hline
    &      &                                       &      & 2/3                          & 9    & 18  \\ \hline
40  & 2    & ${A_2}^2 \bigoplus {A_1}^2$           & 16   & 3/2                          & 2    & 4   \\ \hline
    &      &                                       &      & 4/3                          & 3    & 4   \\ \hline
    &      &                                       &      & 5/6                          & 6    & 8   \\  \hline
    &      &                                       &      & 2/3                          & 9    & 4   \\  \hline
    &      &                                       &      & 1/3                          & 12   & 4   \\  \hline
    &      &                                       &      & 1/6                          & 18   & 4   \\  \hline
41  & 2    & $A_2 \bigoplus {A_1}^4$               & 14   & 0                            & 16   & 1   \\  \hline
    &      &                                       &      & 4/3                          & 3    & 6   \\  \hline
    &      &                                       &      & 1                            & 4    & 12  \\  \hline
    &      &                                       &      & 1/3                          & 12   & 12  \\  \hline
42  & 2    & ${A_1}^6$                             & 12   & 0                            & 16   & 3   \\  \hline
    &      &                                       &      & 2                            & 1    & 4   \\  \hline
    &      &                                       &      & 1                            & 4    & 12  \\  \hline
    &      &                                       &      & 1/2                          & 8    & 16  \\  \hline
43  & 1    & $E_7$                                 & 126  & 2                            & 1    & 2   \\  \hline
    &      &                                       &      & 1/2                          & 56   & 2   \\  \hline
44  & 1    & $A_7$                                 & 56   & 0                            & 70   & 1   \\  \hline
    &      &                                       &      & 2                            & 1    & 2   \\  \hline
    &      &                                       &      & 1/2                          & 28   & 4   \\  \hline
45  & 1    & $A_7$                                 & 56   & 9/8                          & 8    & 2   \\  \hline
    &      &                                       &      & 1/2                          & 28   & 2   \\  \hline
    &      &                                       &      & 1/8                          & 56   & 2   \\  \hline
46  & 1    & $D_7$                                 & 84   & 1                            & 14   & 2   \\  \hline
    &      &                                       &      & 1/4                          & 64   & 2   \\  \hline
47  & 1    & $A_6 \bigoplus A_1$                   & 44   & 8/7                          & 7    & 2   \\ \hline
    &      &                                       &      & 9/14                         & 14   & 2   \\  \hline
    &      &                                       &      & 2/7                          & 35   & 2   \\  \hline
    &      &                                       &      & 1/14                         & 42   & 2   \\  \hline
\end{tabular}
}
\end{table}

\begin{table}[ht]
    \centering
    \footnotesize
    \renewcommand{\arraystretch}{1.0}
    \setlength{\tabcolsep}{10pt}
    \caption{Main table}
    \label{tab8:my_label}
    \scalebox{1.1}{
\begin{tabular}{ccccccc}
\hline
No. & $r$  & T                                     & $v(T)$ & $<P_i,P_i>$ & $m(P_i)$ & ${\mathrm \#} {\mathcal{P}_i}$ \\ \hline
48  & 1    & $D_6 \bigoplus A_1$                   & 62   & 0                            & 64   & 1   \\  \hline
    &      &                                       &      & 2                            & 1    & 2   \\  \hline
    &      &                                       &      & 1/2                          & 32   & 2   \\  \hline
    &      &                                       &      & 1/2                          & 24   & 2   \\  \hline
49  & 1    & $E_6 \bigoplus A_1$                   & 74   & 3/2                          & 2    & 2   \\  \hline
    &      &                                       &      & 2/3                          & 27   & 2   \\  \hline
    &      &                                       &      & 1/6                          & 54   & 2   \\  \hline
50  & 1    & $D_5 \bigoplus A_2$                   & 46   & 4/3                          & 3    & 2   \\  \hline
    &      &                                       &      & 3/4                          & 16   & 2   \\  \hline
    &      &                                       &      & 1/3                          & 30   & 2   \\  \hline
    &      &                                       &      & 1/12                         & 48   & 2   \\  \hline
51  & 1    & $A_5 \bigoplus A_2$                   & 36   & 0                            & 45   & 2   \\  \hline
    &      &                                       &      & 2                            & 1    & 2   \\  \hline
    &      &                                       &      & 1/2                          & 20   & 2   \\  \hline
    &      &                                       &      & 1/2                          & 18   & 4   \\   \hline
52  & 1    & $D_5 \bigoplus {A_1}^2$               & 44   & 0                            & 40   & 1   \\  \hline
    &      &                                       &      & 1                            & 10   & 2   \\  \hline
    &      &                                       &      & 1                            & 4    & 2   \\  \hline
    &      &                                       &      & 1/4                          & 32   & 4   \\ \hline
53  & 1    & $A_5 \bigoplus {A_1}^2$               & 34   & 0                            & 40   & 1   \\ \hline
    &      &                                       &      & 3/2                          & 2    & 2   \\ \hline
    &      &                                       &      & 2/3                          & 15   & 2   \\ \hline
    &      &                                       &      & 2/3                          & 12   & 2   \\ \hline
    &      &                                       &      & 1/6                          & 30   & 2   \\ \hline
    &      &                                       &      & 1/6                          & 24   & 2   \\ \hline
54  & 1    & $D_4 \bigoplus A_3$                   & 36   & 0                            & 48   & 1   \\ \hline
    &      &                                       &      & 1                            & 8    & 2   \\ \hline
    &      &                                       &      & 1                             & 6    & 2   \\ \hline
    &      &                                       &      & 1/4                          & 32   & 4   \\ \hline
\end{tabular}
}
\end{table}

\begin{table}[ht]
    \centering
    \footnotesize
    \renewcommand{\arraystretch}{1.0}
    \setlength{\tabcolsep}{10pt}
    \caption{Main table}
    \label{tab9:my_label}
    \scalebox{1.1}{
\begin{tabular}{ccccccc}
\hline
No. & $r$  & T                                     & $v(T)$ & $<P_i,P_i>$ & $m(P_i)$ & ${\mathrm \#} {\mathcal{P}_i}$ \\ \hline
55  & 1    & $A_4 \bigoplus A_3$                   & 32   & 5/4                          & 4    & 2   \\ \hline
    &      &                                       &      & 4/5                          & 10   & 2   \\ \hline
    &      &                                       &      & 9/20                         & 20   & 2   \\ \hline
    &      &                                       &      & 1/5                          & 30   & 2   \\ \hline
    &      &                                       &      & 1/20                         & 40   & 2   \\ \hline
56  & 1    & $A_4 \bigoplus A_2 \bigoplus A_1$     & 28   & 6/5                          & 5    & 2   \\ \hline
    &      &                                       &      & 5/6                          & 6    & 2   \\ \hline
    &      &                                       &      & 8/15                         & 15   & 2   \\ \hline
    &      &                                       &      & 3/10                         & 20   & 2   \\ \hline
    &      &                                       &      & 2/15                         & 30   & 2   \\ \hline
    &      &                                       &      & 1/30                         & 30   & 2   \\ \hline
57  & 1    & $D_4 \bigoplus {A_1}^3$               & 30   & 0                            & 32   & 3   \\ \hline
    &      &                                       &      & 2                            & 1    & 2   \\ \hline
    &      &                                       &      & 1/2                          & 16   & 6   \\ \hline
    &      &                                       &      & 1/2                          & 8    & 2    \\ \hline
58  & 1    & ${A_3}^2 \bigoplus A_1$               & 26   & 0                            & 36   & 3   \\ \hline
    &      &                                       &      & 2                            & 1    & 2   \\ \hline
    &      &                                       &      & 1/2                          & 12   & 6   \\ \hline
    &      &                                       &      & 1/2                          & 16   & 2   \\ \hline
59  & 1    & $A_3 \bigoplus A_2 \bigoplus {A_1}^2$ & 22   & 0                            & 24   & 1   \\ \hline
    &      &                                       &      & 4/3                          & 3    & 2   \\ \hline
    &      &                                       &      & 3/4                          & 8    & 4   \\ \hline
    &      &                                       &      & 1/12                          & 24   & 4   \\ \hline
    &      &                                       &      & 1/3                          & 18   & 2   \\ \hline
    &      &                                       &      & 1/3                          & 12   & 2   \\ \hline
60  & 1    & $A_3 \bigoplus {A_1}^4$               & 20   & 0                            & 24   & 2   \\ \hline
    &      &                                       &      & 0                             & 16   & 1   \\ \hline
    &      &                                       &      & 1                            & 6    & 2   \\ \hline
    &      &                                       &      & 1                             & 4    & 4   \\ \hline
    &      &                                       &      & 1/4                          & 16   & 8   \\ \hline
\end{tabular}
}
\end{table}

\begin{table}[ht]
    \centering
    \footnotesize
    \renewcommand{\arraystretch}{1.0}
    \setlength{\tabcolsep}{10pt}
    \caption{Main table}
    \label{tab10:my_label}
    \scalebox{1.1}{
\begin{tabular}{ccccccc}
\hline
No. & $r$  & T                                     & $v(T)$ & $<P_i,P_i>$ & $m(P_i)$ & ${\mathrm \#} {\mathcal{P}_i}$ \\ \hline
61  & 1    & ${A_2}^3 \bigoplus A_1$               & 20   & 0                            & 27   & 2   \\ \hline
    &      &                                       &      & 3/2                          & 2    & 2   \\ \hline
    &      &                                       &      & 2/3                          & 9    & 6   \\ \hline
    &      &                                       &      & 1/6                          & 18   & 6   \\ \hline
62  & 0    & $E_8$                                 & 240    & no exist                     & 0   & 0   \\ \hline
63  & 0    & $A_8$                                 & 72     & 0                            & 84   & 2   \\ \hline
64  & 0    & $D_8$                                 & 112    & 0                            & 128   & 1   \\ \hline
65  & 0    & $E_7 \bigoplus A_1$                   & 128    & 0                            & 112   & 1   \\ \hline
66  & 0    & $A_5 \bigoplus A_2 \bigoplus A_1$     & 38     & 0                            & 45   & 2   \\ \hline
    &      &                                       &        & 0                            & 40   & 1   \\ \hline
    &      &                                       &        & 0                            & 36   & 2   \\ \hline
67  & 0    & ${A_4}^2$                             & 40     & 0                            & 50   & 4   \\ \hline
68  & 0    & ${A_2}^4$                             & 24     & 0                            & 27   & 8   \\ \hline
69  & 0    & $E_6 \bigoplus A_2$                   & 78     & 0                            & 81   & 2   \\ \hline
70  & 0    & $A_7 \bigoplus A_1$                   & 58     & 0                            & 70   & 1   \\ \hline
    &      &                                       &        & 0                            & 56   & 2   \\ \hline
71  & 0    & $D_6 \bigoplus {A_1}^2$               & 64     & 0                            & 48   & 2   \\ \hline
    &      &                                       &        & 0                            & 64   & 1   \\ \hline
72  & 0    & $D_5 \bigoplus A_3$                   & 52     & 0                            & 60   & 1   \\ \hline
    &      &                                       &        & 0                            & 64   & 2   \\ \hline
73  & 0    & ${D_4}^2$                             & 48     & 0                            & 64   & 3   \\ \hline
74  & 0    & ${A_3}^2 \bigoplus {A_1}^2$           & 28     & 0                            & 24   & 3  \\ \hline
    &      &                                       &        & 0                            & 32   & 1   \\ \hline
    &      &                                       &        & 0                            & 36   & 3   \\ \hline
\end{tabular}
}
\end{table}

\end{document}